\documentclass{article}

\usepackage{amsmath, amssymb, amsthm}

\usepackage{xparse,graphicx,caption,float}
\usepackage[export]{adjustbox}

\newcommand{\sg}{\sigma}

\newtheorem{theorem}{Theorem}

\title{A note on passing from a quasi-symmetric function 
expansion to a Schur function expansion of a symmetric function}

\author{
Adriano M. Garsia \\
\small Department of Mathematics\\[-0.8ex]
\small University of California, San Diego\\[-0.8ex]
\small La Jolla, CA 92093-0112. USA\\[-0.8ex]
\small \texttt{garsiaadriano@gmail.com}
\and
\and
Jeffrey B. Remmel \\
\small Department of Mathematics\\[-0.8ex]
\small University of California, San Diego\\[-0.8ex]
\small La Jolla, CA 92093-0112. USA\\[-0.8ex]
\small \texttt{jremmel@ucsd.edu}
\and
}

\begin{document}

\maketitle 

\begin{center}
{\it Shortly after the  Egge, Loehr and Warrington paper \cite{ELW} became available Jeff Remmel presented the  contents in his topics course. I happened to be in the audience. During  Remmel's presentation it occurred to me that their result implied that when a symmetric function  has been given an expansion in terms of the Gessel fundamentals indexed by compositions then a schur function expansion can be obtained by replacing each Gessel fundamental by a Schur function  indexed by the same composition. After the lecture I deviced the direct proof of this result given in this paper. Upon reading my write up Jeff discovered an error in my involution  and corrected it. Jeff wrote this paper after we encountered a great deal of scepticism about this interpretation of the Egge, Loehr and Warrington result.  The paper remained in my files for several years. After Jeff's passing I decided that this contribution of Jeff should be recorded. The only addition to Jeff exposition I have inserted is some applications of the Egge, Loehr and Warrington result that have been made under my direction and under the direction of Jeff Remmel. These are listed at the end of this manuscript. They include the work of Emily Sergel  \cite{Sergel}, the work of Dun Qiu \cite{Qiu} and work of Austin Roberts \cite{Roberts}}.
\end{center}
 
\begin{abstract}
Egge, Loehr and Warrington  
 gave in \cite{ELW} a combinatorial formula that permits to convert the expansion of a symmetric function,   homogeneous of degree 
$n$,  in terms of Gessel's fundamental quasisymmetric functions into an expansion in terms of Schur functions. Surprisingly the  Egge, Loehr and Warrington result may be shown to be simply equivalent to replacing the Gessel fundamental by a Schur function indexed by the same composition. In this paper we give a direct proof of the validity of this replacement. This interpretation of the result in 
\cite{ELW} has already been successfully applied to 
Schur positivity problems.     
\end{abstract}

\section{Preliminaries}

We say that a sequence of positive integers $\alpha = 
(\alpha_1, \ldots, \alpha_k)$ is 
a composition of $m$ into $k$ parts 
if $\sum_{i=1}^k \alpha_i =m$.  If, in addition, 
$\alpha_1 \geq \ldots \geq \alpha_k$, then we say that 
$\alpha$ is a partition of $m$. 
We say that a sequence of non-negative integers $\gamma=(\gamma_1, \ldots, 
\gamma_{\ell})$ is 
a weak composition of $m$ into $\ell$ parts 
if $\sum_{i=1}^{\ell}\gamma_i =m$.   Thus the difference 
between compositions and weak compositions is that 0 parts are 
allowed in weak compositions.  We shall write 
$\lambda \vdash m$ to denote that $\lambda$ is a partition of $m$, 
$\alpha \Vdash m$ to denote that $\alpha$ is a composition 
of $m$, and $\gamma \Vdash_w m$ to denote the $\gamma$ is 
a weak composition of $m$. Let $S_n$ denote the symmetric group.

Suppose that  
$\gamma = (\gamma_1, \ldots, \gamma_n)$ is a weak compostion of 
$n$ into $n$ parts. We let $X =(x_1, \ldots, x_n)$ and 
\begin{equation*}
\Delta_{\gamma}(X) =  \det||x_i^{\gamma_j+n-j}|| =  \sum_{\sg  \in S_n} sgn(\sg) \sg(x_1^{\gamma_1+n-1} \cdots 
x_n^{\gamma_n+n-n}). 
\end{equation*}
We let $\displaystyle \Delta(X) = \det||x_i^{n-j}||$  
be the Vandermonde determinant. Then the Schur function $s_{\gamma}(X)$  is defined to be 
\begin{equation}
s_{\gamma}(X) = \frac{\Delta_{\gamma}(X)}{\Delta(X)}.
\end{equation}
It is well known that for any such weak composition $\gamma$,  
either $s_{\gamma}(X) =0$ or there is a 
partition $\lambda$ of $n$ such that 
$s_{\gamma}(X) = \pm s_{\lambda}(X)$.
In fact, there is a well-known straightening relation 
which allows to prove that fact. Namely, 
if $\gamma_{i+1}  > 0$, then 
\begin{equation}\label{straight}
s_{(\gamma_1, \ldots, \gamma_i, \gamma_{i+1}, \ldots, \gamma_n)}(X) = - s_{(\gamma_1, \ldots, \gamma_{i+1}-1, \gamma_{i}+1, \ldots, \gamma_n)}(X).
\end{equation}
See \cite{Macbook}.

Suppose $\alpha = (\alpha_1, \ldots, \alpha_k)$ is 
a composition of $n$ with $k$ parts. We associate 
a subset $S(\alpha)$ of $\{1, \ldots, n-1\}$ with $\alpha$ 
by setting 
$$Set(\alpha) = \{\alpha_1, \alpha_1+\alpha_2, \ldots , \alpha_1 + \cdots + 
\alpha_{k-1}\}.$$
We let $\tilde{\alpha}$ be the weak composition 
of $n$ with $n$ parts by adding a sequence of $n-k$ 0's at the 
end of $\alpha$.  For example, if 
$\alpha =(2,3,2,1)$, then $Set(\alpha) =\{2,5,7\}$ and 
$\tilde{\alpha} = (2,3,2,1,0,0,0,0)$. 
Gessel \cite{Ges} introduced 
a fundamental quasisymmetric function associated 
with each composition $\alpha$ which is defined by  
\begin{equation}
F_{\alpha}(X) = 
\sum_{\stackrel{1 \leq a_1 \leq a_2 \leq \cdots \leq a_n \leq n}{i \in 
Set(\alpha) \rightarrow a_i < a_{i+1}}} x_{a_1} x_{a_2} \cdots x_{a_n}.
\end{equation}
The $F_{\alpha}(X)$'s as $\alpha$ ranges over 
the compostions of $k$ are  a basis for the space of quasisymmetric 
functions  $Q_k(x_1, \ldots,x_n)$ of degree $k$. 

There are 
many examples in the literature where one can give a 
combinatorial description of the coefficients that arise 
in the expansion of important symmetric functions 
in terms 
of the fundamental quasisymmetric functions  where one 
does not have a combinatorial interpretation of the 
coefficients in terms of the Schur functions. For 
example, Haglund, Haiman, and Loehr  \cite{HHL} 
gave a combinatorial description of 
the coefficients that arise in expanding 
the modified Macdonald polynomials 
$\tilde{H}_{\mu}(x_1, \ldots, x_n;q,t)$ as a sum of fundamental 
quasi-symmetric functions. Similarly, Loehr and Warrington gave a combinatorial 
description of the coefficients that arise in 
exanding the plethsym of two Schur functions 
in term of fundamental quasisymmetric functions \cite{LW1}. 
The shuffle conjecture of 
Haglund, Haiman, Loehr, Remmel, and Uylanov \cite{HHLRU} provides 
a  conjectured combinatorial description of the 
expansion of the Frobenius  image of 
the character generating function of the space of diagonal 
harmonics in terms parking functions 
weighted by fundamental quasisymmetric functions \cite{HHLRU}. 
In the last few years, there have been several refinements 
and extension of the shuffle conjecture where we have a similar 
situation, see \cite{HMZ}, \cite{GN}, and \cite{BGLX}. 
In all of these cases, we have no combinatorial description 
of the coefficients that arise in the Schur function expansion of 
these symmetric functions.

In a remarkable and important paper, Egge, Loehr and Warrington  
\cite{ELW} gave a combinatorial description of how 
to start with the expansion of a symmetric function  $P(X)$, 
which is homogeneous of degree 
$n$,  in terms of fundamental quasisymmetric functions 
$$P(X) = \sum_{\alpha \Vdash n} a_{\alpha} F_\alpha(X)
$$
and transform it into an expansion in terms of Schur functions 
$$P(X) = 
\sum_{\lambda \vdash n} b_{\lambda} s_{\lambda}(X).
$$
The purpose of this note is to elucidate a simple 
but important consequence of their result. That is, 
we shall prove the following theorem. 

\begin{theorem}\label{main}
Suppose that $P(X)$ is a symmetric function which is 
homogenous of degree $n$ and 
\begin{equation}\label{Pquasi}
P(X) = \sum_{\alpha \Vdash n} a_{\alpha} F_\alpha(X).
\end{equation}
Then 
\begin{equation}\label{Pschur}
P(X) = \sum_{\alpha \Vdash n} 
a_{\alpha} s_{ {\alpha}}(X).
\end{equation}
\end{theorem}
Thus to obtain the Schur function of $P(X)$, one 
simply has to replace each $F_{\alpha}(X)$ 
by $s_{ {\alpha}}(X)$ and then straighten 
the resulting Schur functions. 

As we shall see the proof of Theorem \ref{main} is 
much simplier than the original proof in \cite{ELW}. 

\newpage
\section{Proof of Theorem 1}

We start with the basic fact that if $\mathcal{A}_n$ is the 
polynomial operator  
$$\mathcal{A}_n = \sum_{\sg \in S_n} sgn(\sg) \sg,$$ 
where for any monomial 
$x_1^{a_1} \cdots x_n^{a_n}$ and  
$\sg = \sg_1 \ldots \sg_n \in S_n$, $\sg(x_1^{a_1} \cdots x_n^{a_n})  = 
x_{\sg_1}^{a_1} \cdots x_{\sg_n}^{a_n}$, then 
for any symmetric function $f(X)$ we have, 
\begin{equation}
f(X) = \frac{1}{\Delta(X)} 
\mathcal{A}_n f(X) x^{\delta_n}
\end{equation}
where $x^{\delta_n} = \prod_{i=1}^n x_i^{n-i}$. 
This is an immediate consequence of the determinantal expansion 
$$\Delta(X) = \sum_{\sg \in S_n} sgn(\sg) \sg(x^{\delta_n})$$ 
and the fact that for any $\sg = \sg_1 \ldots \sg_n \in S_n$, 
$$f(X) = f(x_{\sg_1}, \ldots , x_{\sg_n}).$$

Thus to prove Theorem \ref{main}, we need only 
prove that for each composition 
$\alpha =(\alpha_1, \ldots, \alpha_k)$ of $n$, 
\begin{equation}\label{eq:key}
s_{\tilde{\alpha}}(X) 
=  \frac{1}{\Delta(X)} 
\mathcal{A}_n  F_{\alpha}(X) x^{\delta_n} = 
\sum_{\stackrel{1 \leq a_1 \leq a_2 \cdots \leq a_n \leq a_n}{i \in 
Set(\alpha) \rightarrow a_i < a_{i+1}}} 
\frac{1}{\Delta(X)} 
\mathcal{A}_n  x_{a_1} x_{a_2} \cdots x_{a_n} x^{\delta_n}.
\end{equation}

We consider the following involution $I$ of the 
monomials that appear on the right-hand side of (\ref{eq:key}).
First $I$ has one fixed point, namely, the 
monomial $x_1^{\alpha_1} \cdots x_{k}^{\alpha_k} x^{\delta_n}$.
 Given any other monomial  
$x^{u} = x_{a_1} x_{a_2} \cdots x_{a_n} x^{\delta_n}$ which 
appears on the right-hand side of (\ref{eq:key}), look 
for the  $s =s(u) < k$ such that 
$$x^u = x_1^{\alpha_1} x_2^{\alpha_2} 
 \cdots x_s^{\alpha_s} x_{s+1}^{b_{s+1}} \cdots x_{s+r-1}^{b_{r+s-1}} 
x_{s+r}^{b_{s+r}} \left( \prod_{i=s+r+1}^n x_{i}^{b_{i}}\right) x^{\delta_n}$$ 
where $r = r(u) \geq 2$, $b_{s+1} + \cdots +b_{s+r} = \alpha_{s+1}$, 
and $b_{s+r} > 0$. For example, 
if $\alpha = (2,3,3)$, then 
$(a_1, \ldots, a_8) = (1,1,2,2,2,3,5,5)$, then 
$x^u = x_{a_1} \cdots x_{a_8} = x_1^2 x_2^3 x_3^1 x_4^0 x_5^2 
x_6^{0} x_7^{0} x_8^0$ so that $s(u) =2$, $r(u) =3$, and 
$x_{s+1}^{b_{s+1}} \cdots x_{s+r}^{b_{s+r}} = x_3x_4^0x_5^2$. 

Then we let  
\begin{multline}
I\left(x_1^{\alpha_1} x_2^{\alpha_2} 
 \cdots x_s^{\alpha_s} x_{s+1}^{b_{s+1}} \cdots x_{s+r-1}^{b_{r+s-1}}
x_{s+r}^{b_{s+r}} 
\left( \prod_{i=s+r+1}^n x_i^{b_{i}}\right) x^{\delta_n} \right) =\\
x_1^{\alpha_1} x_2^{\alpha_2} 
 \cdots x_s^{\alpha_s} x_{s+1}^{b_{s+1}} \cdots 
x_{s+r-1}^{b_{r+s}-1}x_{s+r}^{b_{s+r-1}+1} 
 \left( \prod_{i=s+r+1}^n x_{i}^{b_{i}}\right)  x^{\delta_n}.
\end{multline}
That is, $I$ simply replaces the factor 
$x_{r+s-1}^{b_{r+s-1}} x_{r+s}^{b_{r+s}}$ in $x^u$ by 
$x_{s+r-1}^{b_{r+s}-1} x_{s+r}^{b_{s+r-1}+1}$. 
It is easy to check  that if 
$$x^v = x_1^{\alpha_1} x_2^{\alpha_2} 
 \cdots x_s^{\alpha_s} x_{s+1}^{b_{s+1}} \cdots 
x_{s+r-1}^{b_{r+s}-1} x_{s+r}^{b_{s+r-1}+1} 
\left( \prod_{i=s+r+1}^n x_{i}^{b_{i}}\right) x^{\delta_n} = 
x_{c_1} \ldots x_{c_n} x^{\delta_n},$$
then $s(v) = s(u)$, $r(u) =r(v)$, and $(c_1, \ldots, c_n)$ 
is a sequence such that $c_1 \leq \cdots \leq c_n$ 
and $i \in S(\alpha)$ implies $c_i < c_{i+1}$ so that 
$x^v$ also appears on the right-hand side of (\ref{eq:key}). 
In addition, it is easy to see that if 
\begin{eqnarray*}
\frac{1}{\Delta(X)} 
\mathcal{A}_n x^u &=& \frac{1}{\Delta(X)} 
\mathcal{A}_n \left(x_1^{\alpha_1} x_2^{\alpha_2} 
 \cdots x_s^{\alpha_s} x_{s+1}^{b_{s+1}} \cdots x_{s+r-1}^{b_{r+s-1}}
x_{s+r}^{b_{s+r}}
\left(\prod_{t > r} x_{s+r}^{b_{s+r}}\right) x^{\delta_n}\right) \\
&=& s_{(\gamma_1, \ldots, \gamma_{r+s-1},\gamma_{r+s}, \ldots, \gamma_n)}(X), 
\end{eqnarray*}
then 
\begin{eqnarray*}
\frac{1}{\Delta(X)} 
\mathcal{A}_n x^v &=& \frac{1}{\Delta(X)} 
\mathcal{A}_n \left(x_1^{\alpha_1} x_2^{\alpha_2} 
 \cdots x_s^{\alpha_s} x_{s+1}^{b_{s+1}} \cdots x_{s+r-1}^{b_{r+s}-1}
x_{s+r}^{b_{s+r-1}+1} 
\left(\prod_{t > r} x_{s+r}^{b_{s+r}}\right) x^{\delta_n}\right) \\
&=& s_{(\gamma_1, \ldots, \gamma_{r+s}-1,\gamma_{r+s-1}+1, \ldots, 
\gamma_n)}(X)  
\end{eqnarray*}
so that by (\ref{straight}), these two terms cancell each other. 

Thus  $I$ shows that the right-hand side of (\ref{eq:key}) reduces 
to 
$$\frac{1}{\Delta(X)} 
\mathcal{A}_n  x_1^{\alpha_1} \ldots x_k^{\alpha_k}  x^{\delta_n} = 
s_{\tilde{\alpha}}(X)$$
which is what we wanted to prove. 

\section{Some applications }

After discovering the present interpretation of the Egge-Loehr-Warrington result. Some efforts were directed towards identifying the surviving terms  after the replacement of a Gessel fundamental by a compositional indexed  Schur function. 
The first successful use of this kind of the Egge-Loehr-Warrington result was obtained  by Emily Sergel in \cite{Sergel}.
Encouraged by  Sergel's  success Dun Qiu and Jeff Remmel, in a truly  remarkable paper \cite{Qiu},  were able to prove  Schur  positivity  for a wider variety of Rational Parking function modules.     
\vskip .1in

We will next describe a specific example were our attempts  led to a conjecture with measurable success.
Let us recall that in \cite{HHL} Haglund, Haiman and Loehr derive  the Lascoux-Schutzenberger  charge result from their combinatorial proof the Haglund formula.  Since their work  consisted in showing that co-charge came out of the Haglund's $inv_\mu$ statistic it was compelling to see if co-charge could be bypassed altogether. This led to the following computer experimentation.

The point of departure is the identity
$$
\widetilde  H_\mu[X;t]\, =\,   \widetilde H_\mu[X;0,t]
\eqno 3.1
$$
expressing a modified Hall-Littlewood polynomial in terms of the modified Macdonald polynomial. Now, in the present context Haglund's formula may be written in the form
$$
\widetilde H_\mu[X;q,t]\, =\, \sum_{\sigma\in S_n}
t^{maj_\mu(\sigma)}q^{inv_\mu(\sigma)}F_{pides(\sigma)}[X]
\eqno 3.2
$$ 
where the french Ferrer's diagram of $\mu$ is filled 
by $\sigma$ in the reading order, that is  by rows from left to right and from top to bottom. The statistic ``$maj_\mu(\sigma)$'' is simply the sum of the major indexes of the column of $\mu$ read from top to bottom,  ``$inv_\mu(\sigma)$'' counts the number 
of counterclockwise triplets and ``$pides(\sigma)$'' gives the composition of the descent set of the inverse of $\sigma$. Thus 3.1 reduces this identity to
$$
\widetilde H_\mu[X;0,t]\, =\, \sum_{\sigma\in S_n; inv_\mu(\sigma)=0}
t^{maj_\mu(\sigma)}F_{pides(\sigma)}[X]
\eqno 3.3
$$ 
Now in an unpublished algorithm  Loehr and Warrington show how to construct $inv_\mu(\sigma)=0$ fillings. Their algorithm is based  on the fact that, for $k=\l(\mu)$,   
it suffices to choose the decomposition 
$$
T_1+T_2+\cdots +T_k=\{1,2,\ldots ,n\}
\hskip .5 in (\hbox{with $|T_i|=\mu_i$)}
$$
of the  entries of $\sigma$ to be placed in the rows 
of $\mu$. In fact, once the first row of $\mu$ is filled by the elements of $T_1$ in increasing order, then the $\,inv_\mu(\sigma)=0 \,$ condition recursively forces the order in which row $i$ must be filled by the elements of $T_i$.
\newpage

 This given, 3.3 may be rewritten as 
$$
\widetilde H_\mu[X;0,t]\, =\,
\sum_{\sigma:T_1+T_2+\cdots +T_k=[1,n]}
t^{maj_\mu(\sigma)}F_{pides(\sigma)}[X]
$$
and our interpretation of the Egge-Loehr-Warrington result gives
$$
\widetilde H_\mu[X;0,t]\, =\,
\sum_{\sigma:T_1+T_2+\cdots +T_k=[1,n]}
t^{maj_\mu(\sigma)}s_{pides(\sigma)}[X].
\eqno 3.4
$$
Now it is well known that we have three alternatives
$$
s_{pides(\sigma)}[X]\, =\, 
\begin{cases}
0                    &\mbox{$s_{pides_(\sigma)}$  straightens to  $0$}  \\
-s_{\lambda(\sigma)} &\mbox{$s_{pides_(\sigma)}$  straightens to  $-s_{\lambda(\sigma)}$} \\                      
s_{\lambda(\sigma)}  &\mbox{$s_{pides_(\sigma)}$  straightens to  $s_{\lambda(\sigma)}$}       
\end{cases}
$$
The parking functions that produce the first alternative do not contribute to the  sum. Due to the Schur positivity of the left hand side of 3.4, the parking functions that produce the second  alternative must cancel out  with exactly one of the parking functions that produces the third  alternative with exactly the same 
$\, maj_\mu(\sigma)\, $ statistic. The resulting sum is over a subset of the original parking functions. An a priori identification of the left overs would deliver the Schur function expansion of the modified Hall-Littlewood polynomial
$\widetilde H_\mu[X;0,t]$.

This given, what initially felt as a wild guess, was the conjecture that the left overs are the $inv_\mu(\sigma)=0$ fillings 
 that produce the third  alternative and that in addition the Schensted row insertion of $\sigma$ results  in a pair of standard tableaux  of shape $\lambda (\sigma)$.

The resulting computer data revealed the astonishing fact that the ``leftover'' according to this simple  criterion yielded the correct  Schur expansion of $\widetilde H_\mu[X;0,t]$ up to partition of $9$ excluding the partition $[3,3,3]$. But even in that case the Schur expansion was only short one term.
the existence of   this counter example discouraged further experimentations. But given the size of the 
counter  example one could be left with the idea that a suitable $\mu$-variant  of the Schensted algorithm may correctly identify  the leftovers without exceptions.

This particular study of the consequences of Haglund's formula entered a new chapter as a result of  
a poster of Austin Roberts in the Paris FPSAC of 2013. This poster exhibited a similar experiment involving the unrestricted Haglund Formula. The Roberts experiment revealed that the Schur expansion  of $\widetilde H_\mu[X;q,t]$  could be obtained from Schensted correspondence provided $\mu$ did not contain 
$[3,3,3]$ and another partition. This circumstance prompted the first named author to ask Roberts to see if in the case of the modified Hall-Littlewood $\widetilde H_\mu[X;0,t]$ the only obstruction to the use of Schensted  to obtain the Schur expansion was containment of $[3,3,3]$. It turned out Roberts succeeded not only in proving this fact but  also showed in  \cite{Roberts} how the conjectured algorithm   had to be modified to yield the correct answer without exceptions. 
\vskip .1 in

Very recently, we received from Ira Gessel a manuscript with a new proof of the Egge-Loehr-Warrington result obtained by constructing an involution that proves the validity of our replacement for the Schur basis.


\newpage
\begin{thebibliography}{19} 





\bibitem{BGLX} F. Bergeron, A. Garsia, E. Leven, and G. Xin, 
A compositional $(km,kn)$-shuffle conjecture, 
arXiv:1404.4616 (2014). 


\bibitem{ELW} E. Egge, N. Loehr, and G. Warrington, 
From quasisymmetric expansion to Schur expansion via 
a modified inverse Kostka matrix, European J. Combin., 
{\bf 31} (2010) no.8, 2014-2027. 


\bibitem{Ges} I. Gessel, Multipartite $P$-partitions and inner products of skew Schur functions, Contemp. Math., {\bf 34} (1984), 289-301. 


\bibitem{GN} E. Gorsky and A. Negut, Refined knot invariants and 
Hilbert Schemes, arXiv:1304.3328 (2013)


\bibitem{HHL} J.~Haglund, M.~Haiman, and N.~Loehr, A Combinatorial 
Formulas for the Macdonald polynomials. Jour. Amer. Math. Soc. {\bf 18} 
(2005), pp. 735-761. 


\bibitem{HHLRU} M. Haiman, J. Haglund, N. Loehr, J.B. Remmel, and A. Ulyanov, 
A combinatorial formula for the character of the diagonal covariants, 
Duke Mathematical Journal, {\bf 126, no. 2}, (2005), 195-232.

\bibitem{HMZ} J. Haglund, J. Morse, and M. Zabrocki, A compositional 
refinement of the shuffle conjecture specifying touch points of the Dyck path, Canadian J. Math, {\bf 64} (2102), 822-844. 



\bibitem{LW1} N. Loehr and G. Warrington, Quasisymmetric expansions of Schur-function 
plethysms, Proceedings of AMS, {\bf 140} (2012), 1159-1171.

\bibitem{Macbook}
I.~G. Macdonald, \emph{Symmetric functions and {H}all polynomials}, second ed.,
  The Clarendon Press, Oxford University Press, New York, 1995, With
  contributions by A.~Zelevinsky, Oxford Science Publications.

\bibitem{Qiu} Dun Qiu and Jeffrey Remmel, 
Schur Function Expansions and the Rational Shuffle Conjecture,  Proceedingsof the $29^{th}$ Conference on Formal Power Series and Algebraic Combinatorics, London (2017).




\bibitem{Roberts} Austin Roberts, On the Schur expansion of Hall-Littlewood  and related polynomials via Yamanouchi words,
arXiv:1404.1036v3 (2015)

\bibitem{Sergel} Emily Sergel, Two Special
Cases of the Rational Shuffle Conjecture, Proceedingsof the $26^{th}$ Conference on Formal Power Series and Algebraic Combinatorics, Chicago (2014).

\end{thebibliography}
\end{document}